# Near approximations via general ordered topological spaces

M.Abo-Elhamayel
Mathematics Department, Faculty of Science
Mansoura University

**Abstract**
Rough set theory is a new mathematical approach to imperfect knowledge. The notion of rough sets is generalized by using an arbitrary binary relation on attribute values in information systems, instead of the trivial equality relation. The topology induced by binary relations is used to generalize the basic rough set concepts. This paper studies near approximation via general ordered topological approximation spaces which may be viewed as a generalization of the study of near approximation from the topological view. The basic concepts of some increasing (decreasing) near approximations, increasing (decreasing) near boundary regions and increasing (decreasing) near accuracy were introduced and sufficiently illustrated. Moreover, proved results, implications and add examples.

## 1. Introduction

The concept of rough set has many applications in data analysis. Topology [5], one of the most important subjects in mathematics, provides mathematical tools and interesting topics in studying information systems and rough sets [2,7,8,11,12,13]. The purpose of this paper is to put a starting point for the applications of ordered topological spaces into rough set analysis. Rough set theory introduced by pawlak in 1982, is a mathematical tool that supports the uncertainty reasoning. Rough sets generalized by many ways [3,6,9,15]. In this paper, we give a general study of $\alpha, P$ approximations, which studied in [1]. Our results in this paper became the results, which obtained before in case of taking the partially ordered relation as an equal relation.

## 2. Preliminaries

In this section, we give an account for the basic definitions and preliminaries to be used in the paper.

**Definition 2.1[10].** A subset $A$ of $U$, where $(U, \rho)$ is a partially ordered set is said to be increasing (resp. decreasing) if for all $a \in A$ and $x \in U$ such that $a\rho x$ (resp. $x\rho a$) imply $x \in A$.

**Definition 2.2[10].** A triple $(U, \tau, \rho)$ is said to be a topological ordered space, where $(U, \tau)$ is a topological apace and $\rho$ is a partially order relation on $U$.

**Definition 2.3[11].** An information system is a pair $(U, \mathbf{A})$, where $U$ is a non-empty finite set of objects and $\mathbf{A}$ is a non-empty finite set of attributes.

**Definition 2.4[4].** A non-empty set $U$ equipped with a general relation $R$ which generate a topology $\tau_R$ on $U$ and a partially order relation $\rho$ wright as $(U, \tau_R, \rho)$ is said to be general ordered topological approximation space (for short, GOTAS).

**Definition 2.5[4].** Let $(U, \tau_R, \rho)$ be a GOTAS and $A \subseteq U$. We define:

(1) $\underline{R}_{Inc}(A) = A^{\circ Inc}$, $A^{\circ Inc}$ is the greatest increasing open subset of $A$.

(2) $\underline{R}_{Dec}(A) = A^{\circ Dec}$, $A^{\circ Dec}$ is the greatest decreasing open subset of $A$.

(3) $\overline{R}^{Inc}(A) = \overline{A}^{Inc}$, $\overline{A}^{Inc}$ is the smallest increasing closed superset of $A$.

(4) $\overline{R}^{Dec}(A) = \overline{A}^{Dec}$, $\overline{A}^{Dec}$ is the smallest decreasing closed superset of $A$.

(5) $\alpha^{Inc} = \dfrac{card(\underline{R}_{Inc}(A))}{card(\overline{R}^{Inc}(A))}$ (resp. $\alpha^{Dec} = \dfrac{card(\underline{R}_{Dec}(A))}{card(\overline{R}^{Dec}(A))}$) and $\alpha^{Inc}$ (resp. $\alpha^{Dec}$), is $R-$ increasing (resp. decreasing) accuracy.

**Definition 2.6[4].** Let $(U, \tau_R, \rho)$ be a GOTAS and $A \subseteq U$. We define:

(1) $\underline{S}_{Inc}(A) = A \cap \overline{R}^{Inc}(\underline{R}_{Inc}(A))$, $\underline{S}_{Inc}(A)$ is called $R-$inc semi lower.

(2) $\overline{S}^{Inc}(A) = A \cup \underline{R}_{Inc}(\overline{R}^{Inc}(A))$, $\overline{S}^{Inc}(A)$ is called $R-$ inc semi upper.

(3) $\underline{S}_{Dec}(A) = A \cap \overline{R}^{Dec}(\underline{R}_{Dec}(A))$, $\underline{S}_{Dec}(A)$ is called $R-$dec semi lower.

(4) $\overline{S}^{Dec}(A) = A \cup \underline{R}_{Dec}(\overline{R}^{Dec}(A))$, $\overline{S}^{Dec}(A)$ is called $R-$dec semi upper.

$A$ is $R-$ increasing (resp. decreasing) semi exact if $\underline{S}_{Inc}(A) = \overline{S}^{Inc}(A)$ (resp. $\underline{S}_{Dec}(A) = \overline{S}^{Dec}(A)$), otherwise $A$ is $R-$ increasing (resp. decreasing) semi rough.

### 3. New approximations and its properties

In this section, we introduce some definitions and propositions about near approximations, near boundary regions via GOTAS, which are essential for present study.

**Definition 3.1.** Let $(U, \tau_R, \rho)$ be a GOTAS and $A \subseteq U$. We define:

(1) $\underline{\alpha}_{Inc}(A) = A \cap \underline{R}_{Inc}(\overline{R}^{Inc}(\underline{R}_{Inc}(A)))$, $\underline{\alpha}_{Inc}(A)$ is called $R$-increasing $\alpha$ lower.

(2) $\overline{\alpha}^{Inc}(A) = A \cup \overline{R}^{Inc}(\underline{R}_{Inc}(\overline{R}^{Inc}(A)))$, $\overline{\alpha}^{Inc}(A)$ is called $R$-increasing $\alpha$ upper.

(3) $\underline{\alpha}_{Dec}(A) = A \cap \underline{R}_{Dec}(\overline{R}^{Dec}(\underline{R}_{Dec}(A)))$, $\underline{\alpha}_{Dec}(A)$ is called $R$-decreasing $\alpha$ lower.

(4) $\overline{\alpha}^{Dec}(A) = A \cup \overline{R}^{Dec}(\underline{R}_{Dec}(\overline{R}^{Dec}(A)))$, $\overline{\alpha}^{Dec}(A)$ is called $R$-decreasing $\alpha$ upper.

$A$ is $R$-increasing (resp. $R$-decreasing) $\alpha$ exact if $\underline{\alpha}_{Inc}(A) = \overline{\alpha}^{Inc}(A)$ (resp. $\underline{\alpha}_{Dec}(A) = \overline{\alpha}^{Dec}(A)$), otherwise $A$ is $R$-increasing (resp. $R$-decreasing) $\alpha$ rough.

**Proposition 3.2.** Let $(U, \tau_R, \rho)$ be a GOTAS and $A, B \subseteq U$. Then

(1) If $A \subseteq B \to \overline{\alpha}^{Inc}(A) \subseteq \overline{\alpha}^{Inc}(B)$ ( $A \subseteq B \to \overline{\alpha}^{Dec}(A) \subseteq \overline{\alpha}^{Dec}(B)$ ).

(2) $\overline{\alpha}^{Inc}(A \cap B) \subseteq \overline{\alpha}^{Inc}(A) \cap \overline{\alpha}^{Inc}(B)$ ($\overline{\alpha}^{Dec}(A \cap B) \subseteq \overline{\alpha}^{Dec}(A) \cap \overline{\alpha}^{Dec}(B)$).

(3) $\overline{\alpha}^{Inc}(A \cup B) \subseteq \overline{\alpha}^{Inc}(A) \cup \overline{\alpha}^{Inc}(B)$ ($\overline{\alpha}^{Dec}(A \cup B) \subseteq \overline{\alpha}^{Dec}(A) \cup \overline{\alpha}^{Dec}(B)$).

**Proof.**

(1) Omitted.

(2) $\overline{\alpha}^{Inc}(A \cap B) = (A \cap B) \cup \overline{R}^{Inc}(\underline{R}_{Inc}(\overline{R}^{Inc}(A \cap B)))$

$\subseteq (A \cap B) \cup \overline{R}^{Inc}(\underline{R}_{Inc}(\overline{R}^{Inc}(A) \cap \overline{R}^{Inc}(B)))$

$\subseteq (A \cap B) \cup \overline{R}^{Inc}(\underline{R}_{Inc}(\overline{R}^{Inc}(A)) \cap \underline{R}_{Inc}(\overline{R}^{Inc}(B)))$

$\subseteq (A \cap B) \cup \overline{R}^{Inc}(\underline{R}_{Inc}(\overline{R}^{Inc}(A)) \cap \overline{R}^{Inc}(\underline{R}_{Inc}(\overline{R}^{Inc}(B))))$

$\subseteq A \cup \overline{R}^{Inc}(\underline{R}_{Inc}(\overline{R}^{Inc}(A)) \cap B \cup \overline{R}^{Inc}(\underline{R}_{Inc}(\overline{R}^{Inc}(B))))$

$\subseteq \overline{\alpha}^{Inc}(A) \cap \overline{\alpha}^{Inc}(B)$.

(3) $\overline{\alpha}^{Inc}(A \cup B) = (A \cup B) \cup \overline{R}^{Inc}(\underline{R}_{Inc}(\overline{R}^{Inc}(A \cup B)))$

$= (A \cup B) \cup \overline{R}^{Inc}(\underline{R}_{Inc}(\overline{R}^{Inc}(A) \cup \overline{R}^{Inc}(B)))$

$\supseteq (A \cup B) \cup \overline{R}^{Inc}(\underline{R}_{Inc}(\overline{R}^{Inc}(A)) \cup \underline{R}_{Inc}(\overline{R}^{Inc}(B)))$

$$\supseteq (A \cup B) \cup \overline{R}^{Inc}(\underline{R}_{Inc}(\overline{R}^{Inc}(A) \cup \overline{R}^{Inc}(\underline{R}_{Inc}(\overline{R}^{Inc}(B))))$$

$$\supseteq A \cup \overline{R}^{Inc}(\underline{R}_{Inc}(\overline{R}^{Inc}(A) \cup B \cup \overline{R}^{Inc}(\underline{R}_{Inc}(\overline{R}^{Inc}(B))))$$

$$\supseteq \overline{\alpha}^{Inc}(A) \cup \overline{\alpha}^{Inc}(B).$$

One can prove the case between parentheses.

**Proposition 3.3.** Let $(U, \tau_R, \rho)$ be a GOTAS and $A, B \subseteq U$. Then

(1) $A \subseteq B \rightarrow \underline{\alpha}_{Inc}(A) \subseteq \underline{\alpha}_{Inc}(B)$ ( $A \subseteq B \rightarrow \underline{\alpha}_{Dec}(A) \subseteq \underline{\alpha}_{Dec}(B)$ ).

(2) $\underline{\alpha}_{Inc}(A \cap B) \subseteq \underline{\alpha}_{Inc}(A) \cap \underline{\alpha}_{Inc}(B)$ ( $\underline{\alpha}_{Dec}(A \cap B) \subseteq \underline{\alpha}_{Dec}(A) \cap \underline{\alpha}_{Dec}(B)$ ).

(3) $\underline{\alpha}_{Inc}(A \cup B) \supseteq \underline{\alpha}_{Inc}(A) \cup \underline{\alpha}_{Inc}(B)$ ( $\underline{\alpha}_{Dec}(A \cup B) \supseteq \underline{\alpha}_{Dec}(A) \cup \underline{\alpha}_{Dec}(B)$ ).

**Proof.**

(1) Easy.

(2) $\underline{\alpha}_{Inc}(A \cap B) = (A \cap B) \cap \underline{R}_{Inc}(\overline{R}^{Inc}(\underline{R}_{Inc}(A \cap B)))$

$$\subseteq (A \cap B) \cap \underline{R}_{Inc}(\overline{R}^{Inc}(\underline{R}_{Inc}(A) \cap \underline{R}_{Inc}(B)))$$

$$\subseteq (A \cap B) \cap \underline{R}_{Inc}(\overline{R}^{Inc}(\underline{R}_{Inc}(A) \cap \overline{R}^{Inc}(\underline{R}_{Inc}(B))))$$

$$\subseteq (A \cap B) \cap \underline{R}_{Inc}(\overline{R}^{Inc}(\underline{R}_{Inc}(A) \cap \underline{R}_{Inc}(\overline{R}^{Inc}(\underline{R}_{Inc}(B)))))$$

$$\subseteq A \cap \underline{R}_{Inc}(\overline{R}^{Inc}(\underline{R}_{Inc}(A) \cap B \cap \underline{R}_{Inc}(\overline{R}^{Inc}(\underline{R}_{Inc}(B)))))$$

$$\subseteq \underline{\alpha}_{Inc}(A) \cap \underline{\alpha}_{Inc}(B).$$

(3) $\underline{\alpha}_{Inc}(A \cup B) = (A \cup B) \cap \underline{R}_{Inc}(\overline{R}^{Inc}(\underline{R}_{Inc}(A \cup B)))$

$$\supseteq (A \cup B) \cap \underline{R}_{Inc}(\overline{R}^{Inc}(\underline{R}_{Inc}(A) \cup \underline{R}_{Inc}(B)))$$

$$\supseteq (A \cup B) \cap \underline{R}_{Inc}(\overline{R}^{Inc}(\underline{R}_{Inc}(A) \cup \overline{R}^{Inc}(\underline{R}_{Inc}(B))))$$

$$\supseteq (A \cup B) \cap \underline{R}_{Inc}(\overline{R}^{Inc}(\underline{R}_{Inc}(A) \cup \underline{R}_{Inc}(\overline{R}^{Inc}(\underline{R}_{Inc}(B)))))$$

$$\supseteq A \cap \underline{R}_{Inc}(\overline{R}^{Inc}(\underline{R}_{Inc}(A) \cup B \cap \underline{R}_{Inc}(\overline{R}^{Inc}(\underline{R}_{Inc}(B)))))$$

$$\supseteq \underline{\alpha}_{Inc}(A) \cup \underline{\alpha}_{Inc}(B).$$

One can prove the case between parentheses.

**Proposition 3.4.** Let $(U, \tau_R, \rho)$ be a GOTAS and $A, B \subseteq U$. If $A$ is $R-$ increasing (resp. decreasing) exact then $A$ is $\alpha-$ increasing (resp. decreasing) exact.

**Proof.**

Let $A$ be $R-$ increasing exact. Then $\overline{R}^{Inc}(A) = \underline{R}_{Inc}(A)$, $\overline{\alpha}^{Inc}(A) = \overline{R}^{Inc}(A)$, $\underline{\alpha}_{Inc}(A) = \underline{R}_{Inc}(A)$. Therefore $\overline{\alpha}^{Inc}(A) = \underline{\alpha}_{Inc}(A)$.

One can prove the case between parentheses.

**Definition 3.5.** Let $(U, \tau_R, \rho)$ be a GOTAS and $A \subseteq U$. Then

(1) $B_{\alpha Inc}(A) = \overline{\alpha}^{Inc}(A) - \underline{\alpha}_{Inc}(A)$ (resp. $B_{\alpha Dec}(A) = \overline{\alpha}^{Dec}(A) - \underline{\alpha}_{Dec}(A)$), is increasing (resp. decreasing) $\alpha$ boundary region.

(2) $Pos_{\alpha Inc}(A) = \underline{\alpha}_{Inc}(A)$ (resp. $Pos_{\alpha Dec}(A) = \underline{\alpha}_{Dec}(A)$), is increasing (resp. decreasing) $\alpha$ positive region.

(3) $Neg_{\alpha Inc}(A) = U - \overline{\alpha}^{Dec}(A)$ (resp. $Neg_{Dec}(A) = U - \overline{\alpha}^{Inc}(A)$), is increasing (resp. decreasing) $\alpha$ negative region.

**Proposition 3.6.** Let $(U, \tau_R, \rho)$ be a GOTAS and $A, B \subseteq U$. Then

(1) $Neg(A) \supseteq Neg_{\alpha Dec}(A)$ ( $Neg(A) \supseteq Neg_{\alpha Dec}(A)$ ).

(2) $Neg_{\alpha Inc}(A \cup B) \subseteq Neg_{\alpha Inc}(A) \cup Neg_{\alpha Inc}(B)$

$( Neg_{\alpha Dec}(A \cup B) \subseteq Neg_{\alpha Dec}(A) \cup Neg_{\alpha Dec}(B) )$.

(3) $Neg_{\alpha Inc}(A \cap B) \supseteq Neg_{\alpha Inc}(A) \cap Neg_{\alpha Inc}(B)$

$( Neg_{\alpha Dec}(A \cap B) \supseteq Neg_{\alpha Dec}(A) \cap Neg_{\alpha Dec}(B) )$.

**Proof.**

(1) Since $\overline{R}(A) \subseteq \overline{R}^{Dec}(A)$, then $U - \overline{R}(A) \supseteq U - \overline{R}^{Dec}(A)$, therefore $Neg(A) \supseteq Neg_{\alpha Inc}(A)$.

(2) $Neg_{\alpha Inc}(A \cup B) = U - [(A \cup B) \cup \overline{R}^{Dec} \underline{R}_{Dec} \overline{R}^{Dec}(A \cup B)]$

$= U - [(A \cup B) \cup \overline{R}^{Dec} \underline{R}_{Dec}(\overline{R}^{Dec}(A) \cup \overline{R}^{Dec}(B))]$

$\subseteq U - [(A \cup B) \cup \overline{R}^{Dec}(\underline{R}_{Dec} \overline{R}^{Dec}(A) \cup \underline{R}_{Dec} \overline{R}^{Dec}(B))]$

$\subseteq U - [(A \cup B) \cup \overline{R}^{Dec}(\underline{R}_{Dec} \overline{R}^{Dec}(A) \cup \overline{R}^{Dec} \underline{R}_{Dec} \overline{R}^{Dec}(B))]$

$\subseteq U - [A \cup \overline{R}^{Dec}(\underline{R}_{Dec} \overline{R}^{Dec}(A) \cup B \cup \overline{R}^{Dec} \underline{R}_{Dec} \overline{R}^{Dec}(B))]$

$\subseteq U - A \cup \overline{R}^{Dec}(\underline{R}_{Dec} \overline{R}^{Dec}(A) \cap U - B \cup \overline{R}^{Dec} \underline{R}_{Dec} \overline{R}^{Dec}(B))]$

$\subseteq Neg_{\alpha Inc}(A) \cap Neg_{\alpha Inc}(B).$

(3) $Neg_{\alpha Inc}(A \cap B) = U - [(A \cap B) \cup \overline{R}^{Dec} \underline{R}_{Dec} \overline{R}^{Dec}(A \cap B)]$

$= U - [(A \cap B) \cup \overline{R}^{Dec} \underline{R}_{Dec}(\overline{R}^{Dec}(A) \cap \overline{R}^{Dec}(B))]$

$\subseteq U - [(A \cap B) \cup \overline{R}^{Dec}(\underline{R}_{Dec} \overline{R}^{Dec}(A) \cap \underline{R}_{Dec} \overline{R}^{Dec}(B))]$

$\subseteq U - [(A \cap B) \cup \overline{R}^{Dec}(\underline{R}_{Dec} \overline{R}^{Dec}(A) \cap \overline{R}^{Dec} \underline{R}_{Dec} \overline{R}^{Dec}(B))]$

$\subseteq U - [A \cup \overline{R}^{Dec}(\underline{R}_{Dec} \overline{R}^{Dec}(A) \cap B \cup \overline{R}^{Dec} \underline{R}_{Dec} \overline{R}^{Dec}(B))]$

$\subseteq U - A \cup \overline{R}^{Dec}(\underline{R}_{Dec} \overline{R}^{Dec}(A) \cup U - B \cup \overline{R}^{Dec} \underline{R}_{Dec} \overline{R}^{Dec}(B))]$

$\subseteq Neg_{\alpha Inc}(A) \cup Neg_{\alpha Inc}(B).$

One can prove the case between parentheses.

**Definition 3.7.** Let $(U, \tau_R, \rho)$ be a GOTAS and $A \subseteq U$. We define:
(1) $\underline{P}_{Inc}(A) = A \cap \underline{R}_{Inc}(\overline{R}^{Inc}(A))$, $\underline{P}_{Inc}(A)$ is called $R-$increasing Pre lower.
(2) $\overline{P}^{Inc}(A) = A \cup \overline{R}^{Inc}(\underline{R}_{Inc}(A))$, $\overline{P}^{Inc}(A)$ is called $R-$increasing Pre upper.
(3) $\underline{P}_{Dec}(A) = A \cap \underline{R}_{Dec}(\overline{R}^{Dec}(A))$, $\underline{P}_{Dec}(A)$ is called $R-$decreasing Pre lower.
(4) $\overline{P}^{Dec}(A) = A \cup \overline{R}^{Dec}(\underline{R}_{Dec}(A))$, $\overline{P}^{Dec}(A)$ is called $R-$decreasing Pre upper.

$A$ is $R$–increasing (resp. $R$–decreasing) Pre exact if $\underline{P}_{Inc}(A) = \overline{P}^{Inc}(A)$ (resp. $\underline{P}_{Dec}(A) = \overline{P}^{Dec}(A)$), otherwise $A$ is $R$–increasing (resp. $R$–decreasing) Pre rough.

**Proposition 3.8.** Let $(U, \tau_R, \rho)$ be a GOTAS and $A, B \subseteq U$. Then

(1) If $A \subseteq B \rightarrow \overline{P}^{Inc}(A) \subseteq \overline{P}^{Inc}(B)$ ( $A \subseteq B \rightarrow \overline{P}^{Dec}(A) \subseteq \overline{P}^{Dec}(B)$ ).

(2) $\overline{P}^{Inc}(A \cap B) \subseteq \overline{P}^{Inc}(A) \cap \overline{P}^{Inc}(B)$ ($\overline{P}^{Dec}(A \cap B) \subseteq \overline{P}^{Dec}(A) \cap \overline{P}^{Dec}(B)$).

(3) $\overline{P}^{Inc}(A \cup B) \subseteq \overline{P}^{Inc}(A) \cup \overline{P}^{Inc}(B)$ ($\overline{P}^{Dec}(A \cup B) \subseteq \overline{P}^{Dec}(A) \cup \overline{P}^{Dec}(B)$).

**Proof.**

(1) Omitted.

(2) $\overline{P}^{Inc}(A \cap B) = (A \cap B) \cup \overline{R}^{Inc}(\underline{R}_{Inc}(A \cap B))$

$= (A \cap B) \cup \overline{R}^{Inc}(\underline{R}_{Inc}(A) \cap \underline{R}_{Inc}(B)))$

$\subseteq (A \cap B) \cup \overline{R}^{Inc}(\underline{R}_{Inc}(A)) \cap \overline{R}^{Inc}(\underline{R}_{Inc}(B))$

$\subseteq A \cup \overline{R}^{Inc}(\underline{R}_{Inc}(A)) \cap B \cup \overline{R}^{Inc}(\underline{R}_{Inc}(B))$

$\subseteq \overline{P}^{Inc}(A) \cap \overline{P}^{Inc}(B)$.

(3) $\overline{P}^{Inc}(A \cup B) = (A \cup B) \cup \overline{R}^{Inc}(\underline{R}_{Inc}(A \cup B))$

$= (A \cup B) \cup \overline{R}^{Inc}(\underline{R}_{Inc}(A) \cup \underline{R}_{Inc}(B)))$

$\supseteq (A \cup B) \cup \overline{R}^{Inc}(\underline{R}_{Inc}(A)) \cup \overline{R}^{Inc}(\underline{R}_{Inc}(B))$

$\supseteq A \cup \overline{R}^{Inc}(\underline{R}_{Inc}(A)) \cup B \cup \overline{R}^{Inc}(\underline{R}_{Inc}(B))$

$\supseteq \overline{P}^{Inc}(A) \cup \overline{P}^{Inc}(B)$.

One can prove the case between parentheses.

**Proposition 3.9.** Let $(U, \tau_R, \rho)$ be a GOTAS and $A, B \subseteq U$. Then:

(1) $A \subseteq B \rightarrow \underline{P}_{Inc}(A) \subseteq \underline{P}_{Inc}(B)$ ( $A \subseteq B \rightarrow \underline{P}_{Dec}(A) \subseteq \underline{P}_{Dec}(B)$ ).

(2) $\underline{P}_{Inc}(A \cap B) \subseteq \underline{P}_{Inc}(A) \cap \underline{P}_{Inc}(B)$ ($\underline{P}_{Dec}(A \cap B) \subseteq \underline{P}_{Dec}(A) \cap \underline{P}_{Dec}(B)$).

(3) $\underline{P}_{Inc}(A \cup B) \supseteq \underline{P}_{Inc}(A) \cup \underline{P}_{Inc}(B)$ ( $\underline{P}_{Dec}(A \cup B) \supseteq \underline{P}_{Dec}(A) \cup \underline{P}_{Dec}(B)$).

**Proof.**

(1) Easy.

(2) $\underline{P}_{Inc}(A \cap B) = (A \cap B) \cap \underline{R}_{Inc}(\overline{R}^{Inc}(A \cap B))$

$= (A \cap B) \cap \underline{R}_{Inc}(\overline{R}^{Inc}(A) \cap \overline{R}^{Inc}(B)))$

$\subseteq (A \cap B) \cap \underline{R}_{Inc}(\overline{R}^{Inc}(A) \cap (\underline{R}_{Inc}(\overline{R}^{Inc}(B)))$

$\subseteq A \cap \underline{R}_{Inc}(\overline{R}^{Inc}(A) \cap B \cap \underline{R}_{Inc}(\overline{R}^{Inc}(B))$

$\subseteq \underline{P}_{Inc}(A) \cap \underline{P}_{Inc}(B)$.

(3) $\underline{P}_{Inc}(A \cup B) = (A \cup B) \cap \underline{R}_{Inc}(\overline{R}^{Inc}(A \cup B))$

$= (A \cup B) \cap \underline{R}_{Inc}(\overline{R}^{Inc}(A) \cup \overline{R}^{Inc}(B)))$

$\supseteq (A \cup B) \cap \underline{R}_{Inc}(\overline{R}^{Inc}(A) \cup (\underline{R}_{Inc}(\overline{R}^{Inc}(B)))$

$\supseteq A \cap \underline{R}_{Inc}(\overline{R}^{Inc}(A) \cup B \cap \underline{R}_{Inc}(\overline{R}^{Inc}(B))$

$\supseteq \underline{P}_{Inc}(A) \cap \underline{P}_{Inc}(B)$.

**Proposition 3.10.** Let $(U, \tau_R, \rho)$ be a GOTAS and $A, B \subseteq U$. If $A$ is $R$-increasing (resp. decreasing) exact then $A$ is $P$-increasing (resp. decreasing) exact.

**Proof.**

Let $A$ be $R$-increasing exact. Then $\overline{R}^{Inc}(A) = \underline{R}_{Inc}(A)$, $\overline{P}^{Inc}(A) = \overline{R}^{Inc}(A)$, $\underline{P}_{Inc}(A) = \underline{R}_{Inc}(A)$. Therefore $\overline{P}^{Inc}(A) = \underline{P}_{Inc}(A)$.

One can prove the case between parentheses.

**Proposition 3.11.** Let $(U, \tau_R, \rho)$ be a GOTAS and $A, B \subseteq U$. Then we have:

(1) $Neg(A) \supseteq Neg_{PInc}(A)$ ( $Neg(A) \supseteq Neg_{PDec}(A)$ ).

(2) $Neg_{PInc}(A \cup B) \subseteq Neg_{PInc}(A) \cup Neg_{PInc}(B)$

( $Neg_{PDec}(A \cup B) \subseteq Neg_{PDec}(A) \cup Neg_{PDec}(B)$ ).

(3) $Neg_{PInc}(A \cap B) \supseteq Neg_{PInc}(A) \cap Neg_{PInc}(B)$

$(Neg_{PInc}(A \cap B) \supseteq Neg_{PInc}(A) \cap Neg_{PInc}(B))$.

**Proof.**

(1) Since $U - \overline{R}^{Dec}(A) \supseteq U - A \cup \overline{R}^{Dec} \underline{R}_{Dec}(A)$, then $Neg(A) \supseteq Neg_{\alpha Inc}(A)$.

(2) $Neg_{PInc}(A \cup B) = U - [(A \cup B) \cup \overline{R}^{Dec} \underline{R}_{Dec}(A \cup B)]$

$\subseteq U - [(A \cup B) \cup \overline{R}^{Dec}(\underline{R}_{Dec}(A) \cup \underline{R}_{Dec}(B))]$

$\subseteq U - [(A \cup B) \cup \overline{R}^{Dec}(\underline{R}_{Dec}(A) \cup \overline{R}^{Dec} \underline{R}_{Dec}(B))]$

$\subseteq U - [A \cup \overline{R}^{Dec}(\underline{R}_{Dec}(A) \cup B \cup \overline{R}^{Dec} \underline{R}_{Dec}(B))]$

$\subseteq U - A \cup \overline{R}^{Dec}(\underline{R}_{Dec}(A) \cap U - B \cup \overline{R}^{Dec} \underline{R}_{Dec}(B))]$

$\subseteq Neg_{PInc}(A) \cap Neg_{PInc}(B)$.

(3) $Neg_{PInc}(A \cap B) = U - [(A \cap B) \cup \overline{R}^{Dec} \underline{R}_{Dec}(A \cap B)]$

$\supseteq U - [(A \cap B) \cup \overline{R}^{Dec}(\underline{R}_{Dec}(A) \cap \underline{R}_{Dec}(B))]$

$\supseteq U - [(A \cap B) \cup \overline{R}^{Dec}(\underline{R}_{Dec}(A) \cap \overline{R}^{Dec} \underline{R}_{Dec}(B))]$

$\supseteq U - [A \cup \overline{R}^{Dec}(\underline{R}_{Dec}(A) \cap B \cup \overline{R}^{Dec} \underline{R}_{Dec}(B))]$

$\supseteq U - A \cup \overline{R}^{Dec}(\underline{R}_{Dec}(A) \cap U - B \cup \overline{R}^{Dec} \underline{R}_{Dec}(B))]$

$\supseteq U - [\overline{P}^{Dec}(A) \cap \overline{P}^{Dec}(B)]$

$\supseteq Neg_{PInc}(A) \cup Neg_{PInc}(B)$.

**Proposition 3.12.** Let $(U, \tau_R, \rho)$ be a GOTAS and $A \subseteq U$. Then

$\underline{R}_{Inc}(A) \subseteq \underline{\alpha}_{Inc}(A) \subseteq \underline{S}_{Inc}(A)$ $(\underline{R}_{Dec}(A) \subseteq \underline{\alpha}_{Dec}(A) \subseteq \underline{S}_{Dec}(A))$.

**Proof.**

Let $x \in \underline{R}_{Inc}(A)$. Then $x \in \overline{R}^{Inc}(\underline{R}_{Inc}(A))$        (i)

Now, we have $x \in A$ and $x \in \underline{R}_{Inc}(\overline{R}^{Inc}(\underline{R}_{Inc}(A)))$. Then $x \in A \cap \underline{R}_{Inc}(\overline{R}^{Inc}(\underline{R}_{Inc}(A)))$, therefore $x \in \underline{\alpha}_{Inc}(A)$. Hence $\underline{R}_{Inc}(A) \subseteq \underline{\alpha}_{Inc}(A)$.     (1)

Since $x \in A \cap \overline{R}^{Inc}(\underline{R}_{Inc}(A))$, then
$$x \in \underline{S}_{Inc}(A) \quad (2)$$
From (1) and (2), we have $\underline{R}_{Inc}(A) \subseteq \underline{\alpha}_{Inc}(A) \subseteq \underline{S}_{Inc}(A)$.

**Proposition 3.13.** Let $(U, \tau_R, \rho)$ be a GOTAS and $A \subseteq U$. Then

$\underline{\alpha}_{Inc}(A) \subseteq \underline{P}_{Inc}(A)$ ($\underline{\alpha}_{Dec}(A) \subseteq \underline{P}_{Dec}(A)$).

**Proof.**

Since $x \in \underline{R}_{Inc}(\overline{R}^{Inc}(A))$, then $x \in \underline{\alpha}_{Inc}(A)$, and then $x \in A \cap \underline{R}_{Inc}(\overline{R}^{Inc}(\underline{R}_{Inc}(A)))$, therefore $x \in A$ and $x \in \underline{R}_{Inc}(\overline{R}^{Inc}(\underline{R}_{Inc}(A))) \subseteq \underline{R}_{Inc}(\overline{R}^{Inc}(\overline{R}^{Inc}(A)))$. Thus $x \in \underline{R}_{Inc}(\overline{R}^{Inc}(A))$, and thus $x \in A \cap \underline{R}_{Inc}(\overline{R}^{Inc}(A))$. Hence $x \in \underline{P}_{Inc}(A)$.

**Proposition 3.14.** Let $(U, \tau_R, \rho)$ be a GOTAS and $A \subseteq U$. Then $\overline{S}^{Inc}(A) \subseteq \overline{\alpha}^{Inc}(A) \subseteq \overline{R}^{Inc}(A)$ ($\overline{S}^{Dec}(A) \subseteq \overline{\alpha}^{Dec}(A) \subseteq \overline{R}^{Dec}(A)$).

**Proof.**

Let $x \in \overline{S}^{Inc}(A)$, then $x \in A$ or $x \in \underline{R}_{Inc}(\overline{R}^{Inc}(A))$. Thus $x \in A \cup \overline{R}^{Inc}(\underline{R}_{Inc}(\overline{R}^{Inc}(A)))$. Hence
$$x \in \overline{\alpha}^{Inc}(A) \quad (1)$$
Since $x \in A \cup \overline{R}^{Inc}(\underline{R}_{Inc}(\overline{R}^{Inc}(A)))$, then $x \in A \cup \overline{R}^{Inc}(\overline{R}^{Inc}(A))$, therefore $x \in A \cup \overline{R}^{Inc}(A)$. Thus
$$x \in \overline{R}^{Inc}(A) \quad (2)$$
From (1) and (2), we have $\overline{S}^{Inc}(A) \subseteq \overline{\alpha}^{Inc}(A) \subseteq \overline{R}^{Inc}(A)$.

**Definition 3.15.** Let $(U, \tau_R, \rho)$ be a GOTAS and $A \subseteq U$. Then:
(1)    $B_{PInc}(A) = \overline{P}^{Inc}(A) - \underline{P}_{Inc}(A)$ (resp.    $B_{PDec}(A) = \overline{P}^{Dec}(A) - \underline{P}_{Dec}(A)$),    is increasing (resp. decreasing) near boundary region.
(2)    $Pos_{PInc}(A) = \underline{P}_{Inc}(A)$ (resp.    $Pos_{PDec}(A) = \underline{P}_{Dec}(A)$), is increasing (resp.

decreasing) near positive region.

(3) $Neg_{PInc}(A) = U - \overline{P}^{Dec}(A)$ (resp. $Neg_{Dec}(A) = U - \overline{P}^{Inc}(A)$), is increasing (resp. decreasing) near negative region.

**Definition 3.16.** Let $(U, \tau_R, \rho)$ be a GOTAS and $A$ non-empty finite subset of $U$. Then the increasing (decreasing) near accuracy of a finite non-empty subset $A$ of $U$ is given by:

$$\eta_{jInc}(A) = \frac{|\underline{j}_{Inc}(A)|}{|\overline{j}^{Inc}(A)|}, \quad j \in \{\alpha, P\}.$$

**Proposition 3.17.** Let $(U, \tau_R, \rho)$ be a GOTAS and $A$ non-empty finite subset of $U$. Then $\eta_{Inc}(A) \leq \eta_{jInc}(A)$ ($\eta_{Dec}(A) \leq \eta_{jDec}(A)$), for all $j \in \{\alpha, P\}$, where $\eta_{Inc}(A) = \frac{|\underline{R}_{Inc}(A)|}{|\overline{R}^{Inc}(A)|}$ and $\eta_{Dec}(A) = \frac{|\underline{R}_{Dec}(A)|}{|\overline{R}^{Dec}(A)|}$.

**Proof.** Omitted.

**Example 3.18.** Let $U = \{a, b, c, d\}$, $U/R = \{\{a\}, \{a, b\}, \{c, d\}\}$, $\tau_R = \{U, \phi, \{a, b\}, \{c, d\}, \{a\}, \{a, d, c\}\}$, $\tau_R^C = \{U, \phi, \{c, d\}, \{a, b\}, \{b, c, d\}, \{b\}\}$ and $\rho = \{(a, a), (b, b), (c, c), (d, d), (a, b), (b, d), (a, d), (a, c), (c, d)\}$.

For $A = \{a, c\}$, we have:

$\underline{R}_{Dec}(A) = \{a\}$, $\overline{R}^{Dec}(\underline{R}_{Dec}(A)) = \{a, b\}$, $\overline{R}^{Dec}(A) = U$, $\underline{R}_{Dec}(\overline{R}^{Dec}(A)) = U$.

$\underline{P}_{Dec}(A) = A \cap U = A$, $\overline{P}^{Dec}(A) = \{a, b, c\}$, $B_{PDec}(A) = \{b\}$, $Neg_{Inc} = \{d\}$.

$\underline{\alpha}_{Dec}(A) = A \cap \{a, b\} = \{a\}$, $\overline{\alpha}^{Dec}(A) = U$, $B_{\alpha Dec}(A) = \{b, c, d\}$, $Neg_{\alpha Inc} = \phi$.

**Proposition 3.19.** Let $(U, \tau_R, \rho)$ be a GOTAS and $A \subseteq U$. Then we have $B_{SInc}(A) \subseteq B_{\alpha Inc}(A) \subseteq B_{Inc}(A)$ ($B_{SDec}(A) \subseteq B_{\alpha Dec}(A) \subseteq B_{Dec}(A)$).

**Proof.** Omitted.

## 4. Conclusion

As a step, which is rich in results up till now to generalize the generalized approximation spaces, it was the study of GOTAS which is a generalization of the study of OTAS, GAS and AS. Every GOTAS can be regarded as an OTAS if $R$ is an equivalence relation and OTAS can be regarded as an AS if $\rho$ is the equal relation. In addition, every GOTAS can be regarded as GAS if $\rho$ is the equal relation and GAS can be regarded as AS if $R$ is an equivalence relation.

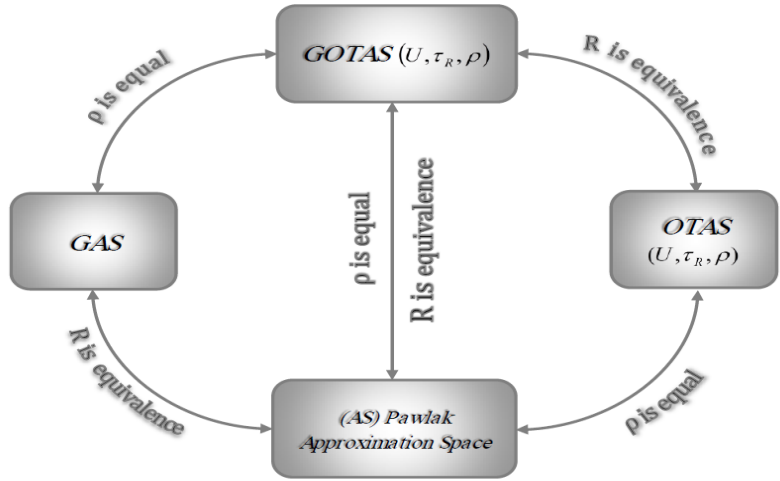